\documentclass[a4paper,12pt, oneside]{amsart}
\usepackage{geometry,graphicx,wrapfig,color}
\usepackage{amssymb}
\usepackage{url}
\def\PL{\mathrm{PL}}
\def\op{\mathrm{op}}
\def\R{\mathbb{R}}
\def\Ord{\mathop{\mathrm{Ord}}}
\def\id{\mathrm{id}}
\def\supp{\mathop{\mathrm{supp}}}
\def\ra{\rangle}
\def\la{\langle}
\def\St{\mathop{\mathrm{Star}}}
\def\Lk{\mathop{\mathrm{Link}}}
\newcommand{\xar}[1]{\ensuremath{\xrightarrow{#1}}}
\newcommand{\wt}[1]{\widetilde{#1}}
\newcommand{\mr}[1]{\mathrm{#1}}
\newcommand{\mb}[1]{\mathbf{#1}}
\newcommand{\bs}[1]{\boldsymbol{#1}}
\newcommand{\bg}[1]{\boldsymbol{\mathfrak{#1}}}
\newcommand{\ms}[1]{\mathsf{#1}}
\newcommand{\us}[2]{\underset{#1}{#2}}
\newcommand{\mc}[1]{\mathcal{#1}}
\newtheorem{thm}{Theorem}
\title{Tangent bundle and Gauss functor
of a combinatorial manifold.}
\thanks{Supported by RFBR grant 05-01-00899 and
 S.S. grant 4329.2006.1.}
\begin{document}

\vspace{1cm}
\maketitle
\begin{center}{\textbf{\large{Nikolai Mn\"ev}}}\\
\vspace{.5cm}
\small PDMI RAS, \\ Fontanka 27,
\\Saint-Petersburg, Russia \\
\texttt{mnev@pdmi.ras.ru} \\
\end{center}\mbox{}
\\
\\
This note is as  an extended talk at the workshop
``Combinatorial Geometries and Applications: Oriented Matroids and Matroids''
Marseille-Luminy, November 7-11, 2005.\\
\\
\subsection{Introduction}
To any poset $\ms P$ we associate
 a poset map $\mr E\ms P\xar{\mr T\ms P} \ms P $ accomplished by two sections
$\ms P\xar{s_0,s_\infty}\mr E\ms P$ and a functor $\ms P\xar{
\mb{G}} \mb{Posets}$. Functor $\mb G$ is related to the map $\mr
T\ms P$ in such a way that $\mr E\ms P\xar{\mr T\ms P} \ms P$ coincides
with $\mathop{\mr{Hocolim}} \mb  G \xar{\Pi} \ms P$.

In \cite{Mnev:2006} it is shown that in the case when $\ms P$ is a
``strict  abstract  manifold'' (the examples are combinatorial
manifolds and boundary complexes of arbitrary convex polytopes) the
following happens. The order complex of $\mr E\ms P$ is a
combinatorial manifold. The map $B\mr E\ms P\xar{B\mr T\ms P} B\ms
P$ together with two sections $Bs_0, Bs_\infty$ is a Kuiper-Lashoff
$(S^n,0,\infty)$ model for the tangent PL-bundle of $B\ms P$. The
image of  functor  $\mb G$ naturally lives in a nice category
$\bg R_n$. The discrete category $\bg R_n$ can be viewed  as a
discrete replacement of the structure group for the discrete
replacements of $\PL$ fiber bundles with fiber $\R^n$. The main
statement it that  this discretization is exact: $B\bg R_n \approx
B\PL_n$ and $B\ms P\xar{B\mb G} B\bg R_n$ is a homotopy model of PL
Gauss map for the PL manifold $B\ms P$. Here $B*$ is the classifying
space functor, which coincides in the case of posets with the
geometric realization of the order complex, in the case of
categories -- with Milnor geometric realization of the nerve and in the
case of simplicial group with the classifying space for
principal bundles. These results can be viewed as a purely
combinatorial variation of some constructions from
 \cite{Le}, \cite{hatcher:1975}, \cite{steinberger:1986}.
 Also this proves the
$\PL$ double for  MacPherson's conjecture on modeling the real
Grassmanian $B\mr O_n$ by the poset of
oriented matroids \cite{MacPherson:1978},\cite{MacPherson:1991}
(see \cite{MnevZigler:1993}). Here we formulate these results
and briefly discuss their proofs from \cite{Mnev:2006}.
\subsection{Tangent bundle and Gauss map of a $\PL$ manifold} \label{ss02}
Milnor in \cite{Milnor:1961} defined the notion of $n$-dimensional
$\PL$ microbundle, the simplicial structure group $\PL_n$ of
microbundles and the classifying space $B\PL_n$. This theory creates
canonical one-to-one correspondence between isomorphism classes of
$n$-dimensional $\PL$ microbundles on a polyhedron $K$ and homotopy
classes of maps from $K$ to $B\PL_n$. Milnor also defined the notion
of tangent microbundle of a $\PL$ manifold $M^n$. A map
$M^n\xar{G}B\PL_n$ representing tangent microbundle of $M^n$ is
called \emph{Gauss map} of $M^n$ and $BPL_n$ is a
 \emph{$\PL$ Grassmanian}. The space $B\PL_n$ and the Gauss map are defined up to homotopy.
In \cite{KLI} and \cite{KLII} Kuiper and Lashof developed the theory
of models for piecewise-linear $\R^n$-bundles. In particular they
established a canonical one-to-one correspondence between isomorphism
classes of  $n$-dimensional $\PL$ microbundles and
$(S^n,\infty)[(S^n,0,\infty)]$ fiber bundles. A piecewise-linear
$(S^n,\infty)$ fiber bundle is a $\PL$ fiber bundle with fiber $S^n$
and a section labeled by $\infty$. A piecewise-linear
$(S^n,0,\infty)$ fiber bundle is a $\PL$ fiber bundle with fiber
$S^n$ and two sections labeled by $0$ and $\infty$. This sections
should have no points in common.

\subsection{Tangent bundle and Gauss functor of a poset} \label{ss03}
Here we introduce a very general and  useless in its full generality
construction.

Let $\ms P$ be a poset, let $p\in \ms P$. We introduce a notation
$\St p$ for the subposet of $\ms P$ which is a union of all
principal ideals containing $p$. Formally:
$$\St p = \{x \in \ms P | \exists y: p\leq y, x\leq y \}$$
We also introduce a notation
$$\Lk p = \{ x \in \St p | p \not \leq x \}$$

Define a new poset $\mr{E}\ms P$ as follows. Denote by $D\ms P$ a
subset of $\ms P\times \ms P$, formed by all  pairs  $(x,y)$ such
that $\exists z\in \ms P: z\geq x, z\geq y$. Pick a new element
$\infty \not\in \ms P$. Set
$$\mr E\ms P=D\ms P\cup \ms P\times\{\infty\}\subset \ms P\times \ms P\cup \{\infty\}$$
Now we define a partial order on $\mr E\ms P$. Set
$$(x_1,y_1)\leq(x_2,y_2) \Leftrightarrow (x_1\us{\ms P}\leq x_2) \wedge
\begin{cases} y_1 \us{\ms P}{\leq} y_2 & \text{ if } y_1,y_2 \in \ms P \\
              y_1 \in \Lk_\ms P x_2 & \text{ if } y_1 \in \ms P, y_2=\infty \\
              y_1=\infty, y_2=\infty
 \end{cases} $$
Set by $\mr E\ms P\xar{\mr T\ms P} \ms P$ the projection on first
argument. We fix two sections of the poset map $\mr T\ms P$: the
diagonal
 $$\ms P\xar{s_0} \mr E\ms P, s_0(x)=(x,x)\in D\ms P\subset \mr E\ms P $$
 and the ``section at infinity'':
 $$\ms P\xar{s_\infty} \mr E\ms P, s_\infty(x)=(x,\infty)$$
We call the defined above set of data $\la \mr T\ms P, s_0, s_\infty \ra$ the \emph{tangent bundle
of a poset $\ms P$}.\\

To any element $x\in \ms P$ there corresponds a subposet $\mb G_x=(\mr
T\ms P)^{-1}(x) \subset \mr E\ms P$ with induced order. We can
describe the structure of $\mb G_x$:
$$\mb G_x=\{(x,y)| y \in \mr{Star}_\ms P\, x\}\cup
\{(x,\infty)\}$$
$$(x,y_1) \us{\mb G_x}{\leq} (x,y_2) \Leftrightarrow \begin{cases}
y_1 \us{\ms P}{\leq}y_2 & \text{ if } y_1,y_2 \in \ms P \\
 y_1 \in \Lk_\ms P x & \text{ if }      y_2=\infty \\
 y_1 = \infty, y_2 = \infty
\end{cases} $$
We will call
the special diagonal element $(x,x)\in\mb G_x$ the
$0$-element.

To any pair $x_1 \leq x_2$  of comparable elements in $\ms P$ we
associate a poset map
$$\mb G_{x_1}\xar{\mb G_{x_1 \leq x_2}}\mb G_{x_2}$$
defined as follows
$$\mb G_{x_1 \leq x_2}(x_1,y)=
\begin{cases}
(x_2,y) & \text{ if } y\in \St x_2\\
(x_2, \infty) &\text{ if } y\not\in \St x_2\\
(x_2,\infty) &\text{ if } y=\infty
\end{cases}$$
Consider the category $\mb{Posets}^{0,\infty}$ of all posets having
two elements specially labelled. One is element is labelled by $0$ and
another by $\infty$. The morphisms of $\mb{Posets}^{0,\infty}$ are
the poset maps preserving labelled elements. The poset maps $\mb
G_{x_1 \leq x_2}$ preserve the elements marked by $0$ (diagonal)
and by $\infty$. So we can regard $\mb G$ as a functor $\ms P
\xar{\mb G} \mb{Posets}^{0,\infty}$.
We will call $\mb G$   \emph{ the Gauss functor of a poset $\ms P$}. \\

The tangent bundle $\la \mr T\ms P, s_0,s_\infty \ra$ can be
identified with certain classical construction associated with $\mb G$.
The construction is known by the names of ``categorial homotopy
colimit'', or ``Grothendieck construction'' \cite{GJ} or ``double
bar construction'' \cite{May:1975} Probably it's first indication is
in Whitehead's construction of the cone of simplicial map
\cite{Whitehead:1939}. In our situation the construction looks
as follows. Let $\ms P$ be a poset and let $\mb F$ be any functor
$\ms P\xar{\mb F} \mb{Posets}$. With functor $\mb F$  we associate a
new poset $\mr{Hocolim}\, \mb F$ and a poset map $\mr{Hocolim}\, \mb
F\xar{\Pi} \ms P$. Put
$$\mr{Hocolim}\, \mb F =\{(x,y)| x \in \ms P, y \in \mb F_x\}$$
and define $(x_0,y_0)\leq (x_1,y_1) $ iff $x_0 \us{\ms P}{\leq} x_1
$ and $\mb F_{x_0 \us{\ms P}{\leq} x_1}(y_0)\us{\mb F_{x_1}}\leq
y_1$ The projection $\mr{Hocolim}\, \mb F\xar{\Pi} \ms P$ is a
projection on the first argument. With a functor $\ms P\xar{\mb H}
\mb{Posets}^{0,\infty}$ we associate three functors to $\ms P\xar{}
\mb{Posets}$. One is the composite of $\mb H$ and forgetful functor
(we denote it by $\wt{\mb{H}}$). The two others are constant
functors $\bs 0, \bs \infty$,
 sending
entire $\ms P$ to $0$ and entire $\ms P$ to $\infty$. The triad
$$(\mr{Hocolim}\, \wt{\mb{H}}, \mr{Hocolim}\, \bs 0, \mr{Hocolim}\, \bs \infty)$$
is exactly $\mr{Hocolim}\, \wt{\mb{H}}$ together with the graphs of
two sections $s_0$ and $s_\infty$ of  projection $\Pi$. So, the
$\mr{Hocolim}$ of a functor with values in $\mb{Posets}^{0,\infty}$
is naturally equipped with $0$ and $\infty$ sections of projection
$\Pi$

The only essential property of our general construction of the
tangent bundle and Gauss functor of a poset $\ms P$ is that the
tangent bundle $\la \mr E\ms P\xar{\mr T\ms P} \ms P, s_0, s_\infty
\ra$ coincides with $\la \mr{Hocolim}\, \wt{\mb{G}} \xar {\Pi} \ms
P, s_0,s_\infty\ra$.

\subsection{Abstract manifolds} \label{ss04}

Here we define a slight generalization of Alexander's combinatorial
manifold.  For the reference on combinatorial topology of  posets
and ball complexes one may use \cite{Bj}.

A \emph{$\PL$-ball complex} is a pair $(X,  U) $, where $X$ is a
compact Euclidean polyhedron and $ U$ is a covering of $X$ by closed
$\PL$-balls such that the following axioms are
satisfied:\\
{\bf plbc1}: the relative interiors of balls from $ U$ form a partition of $X$.\\
{\bf plbc2}: The boundary of each ball from $ U$ is a union of balls from $ U$.\\
A $\PL$-ball complex is defined up to $\PL$-homeomorphism only by
the combinatorics of adjunctions of its balls. Let $\mc D$ be a
$\PL$-ball complex. Consider the poset $\ms P=\mb P(\mc D)$ of all
its balls ordered by inclusion. Than $\mc D$ is cellular complex
$\PL$-homeomorphic to the complex $(B \ms P , \{B \ms P_{\leq
p}\}_{p\in \ms P})$. Here $\ms P_{\leq p}$ is the principal ideal.
For a poset $\ms Q$ we denote by $B \ms Q$ the geometric realization
of the order complex of $\ms Q$.

This makes possible to define  an abstract $\PL$-ball complex: a
finite poset $\ms P$ is called by \emph{abstract $\PL$-ball complex}
if for any $p\in \ms P$ the polyhedron $B \ms P_{<p}$ is a $\PL$
 sphere.
If $\ms P$ is an abstract ball complex than $(B \ms P , \{B \ms
P_{\leq p}\}_{p\in \ms P})$ is a $\PL$-ball complex
\cite{Bj84}.\footnote{The classical combinatorial characterization
of ball complexes is in topological category, not in $\PL$ category, but for
the $\PL$ case  the theory works without any changes.}  We will call The principal
ideals of an abstract ball complex $\ms P$ by its
\textit{balls}. A poset $\ms P$ is \emph{pure} if all maximal chains
have the same  length. A pure poset $\ms P$ is an {\em abstract
manifold} if both $\ms P$ a $\ms P^\op$ are  abstract $\PL$-ball
complexes. Here $\ms P^\op$ is $\ms P$ with the opposite order. If
$\ms P$ is an  abstract manifold then the simplicial complex $\Ord
\ms P$ is a combinatorial manifold in the classical meaning of
Alexander.

Let $\ms P_0$ $\ms P_1$ be the abstract manifolds. Call a poset map
$\ms P_0\xar{\xi}\ms P_1$ an {\em aggregation morphism} if for any
rank $k$ ball $O$ of $\ms P$ the polyhedron $B\xi^{-1}(O)$ is a
$k$-dimensional $\PL$-ball.

\begin{figure} \includegraphics[scale=.33]{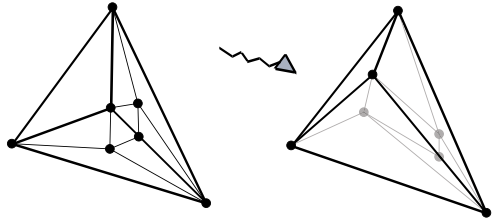}
 \caption{\label{pic1} Abstract aggregation. }\end{figure}

\begin{figure}\includegraphics[scale=.33]{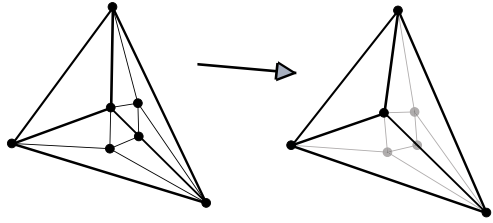}
\caption{\label{pic2}Geometric aggregation.} \end{figure}

An  aggregation morphism $\xi$ \label{aggr} can be  realized up to
isomorphism as a geometric aggregation of $\PL$-ball complex
structures on $B \ms P_0$. Consider the ball complex structure
$(B\ms P_0, \{B {\ms P_0}_{\leq p}\}_{p\in \ms P_0})$ on $B \ms
P_0$. Then after glueing  together all the balls which are sent the same
ball of $\ms P_1$ by morphism $\xi$  we will get a geometric
 representation
of $\ms P_1$ up to an isomorphism. Figure \ref{pic1}
illustrates abstract aggregation and  Figure \ref{pic2} -- geometric
aggregation.
The composition of aggregation morphisms is an aggregation morphism.

Let us call \emph{abstract $n$-sphere} an abstract manifold $\ms
P$ such that polyhedron $B \ms P$ is $\PL$-homeo\-mor\-phic to $S^n$.
Consider the category $\bg R_n$ where the objects are abstract
$n$-spheres having one element of
 maximal rank specially labelled by ``$\infty$''.
 Morphisms of $\bg R_n$ are the
aggregations sending the $\infty$-element  to the $\infty$-element.
So, $\bg R_n$ is a subcategory of the category $\mb{Posets}^\infty$
of all posets with  fixed element labelled by ``$\infty$''.

We need one more definition. An abstract  $n$-dimensional manifold
$\ms P$ is \emph{strict} if for any $x \in \ms P$ the pair of
polyhedra $(B\St x, B \Lk x)$ is $\PL$-homeomorphic to the pair
$(D^n, S^{n-1})$. The examples of strict abstract manifolds are the
simplex posets of combinatorial manifolds and face posets of convex
polytopes.

\subsection{Tangent bundle and Gauss functor of a strict abstract manifold}
Now we can formulate our theorem about the tangent bundle of a poset
in the case of strict abstract manifold. The general message is that
in the  case of strict abstract manifolds (\S\ref{ss04}) the
abstract construction (\S\ref{ss03}) represents the classical
geometrical one (\S\ref{ss02}).

Let $\ms M^n$ be a strict abstract $n$-dimensional manifold. Denote by
$\ms M^n\xar{\mb G^\infty} \mb{Posets^\infty}$ the composition of the
Gauss functor $\mb G$ and forgetful functor $ \mb{Posets}^{0,\infty}
\xar{} \mb{Posets}^{\infty}$.  We will
also call the functor $\mb G^\infty$ ``Gauss functor''.

\begin{thm} \mbox{} \label{th1}\\
1. The image of the Gauss functor $\ms M^n\xar{\mb G^\infty}
\mb{Posets^\infty}$ belongs to
$\bg R_n$ \\
2. The order complex of $\mr E\ms M^n$ is a combinatorial manifold \\
3. The $\PL$-map $B\mr E\ms M^n \xar{B\mr T\ms M^n} B\ms M^n $
together with two sections $Bs_0, Bs_\infty$ is
a Kuiper-Lashof $(S^n, 0,\infty)$ model of the Milnor tangent bundle of $B\ms M^n$.\\
4. The space $B\bg R_n$ is homotopy equivalent to the space $B\PL_n$. \\
5. The map $B\ms M^n\xar{B\mb G^\infty } B\bg R_n$ is a homotopy
model of the $\PL$ Gauss map of $B\ms M^n$.
\end{thm}
If one verifies the statement 1 then the statements 2 and 3  follow
from Kuiper-Lashof theory, simple properties of $\mr{Hocolim}$
construction and Alexander's trick. The statement 5 is a summary of
1-4 modulo standard abstract nonsense. The real thing to prove is the
statement 4 (\cite[Theorem C]{Mnev:2006}).

\begin{figure} $$\includegraphics[scale=0.5]{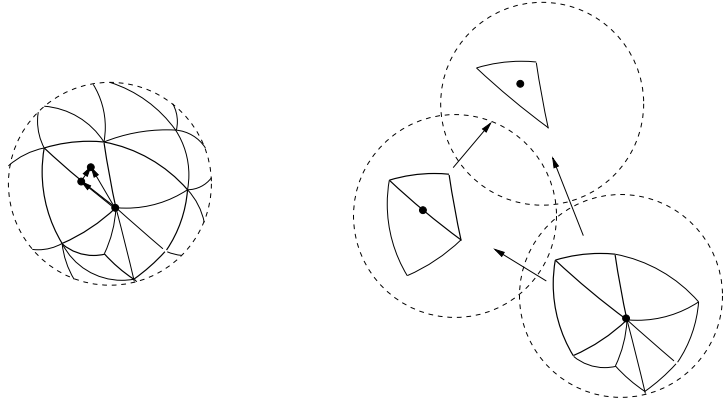}  $$
\caption{\label{pic3} The Gauss functor of a combinatorial sphere}
\end{figure}

\begin{figure}
$$\includegraphics{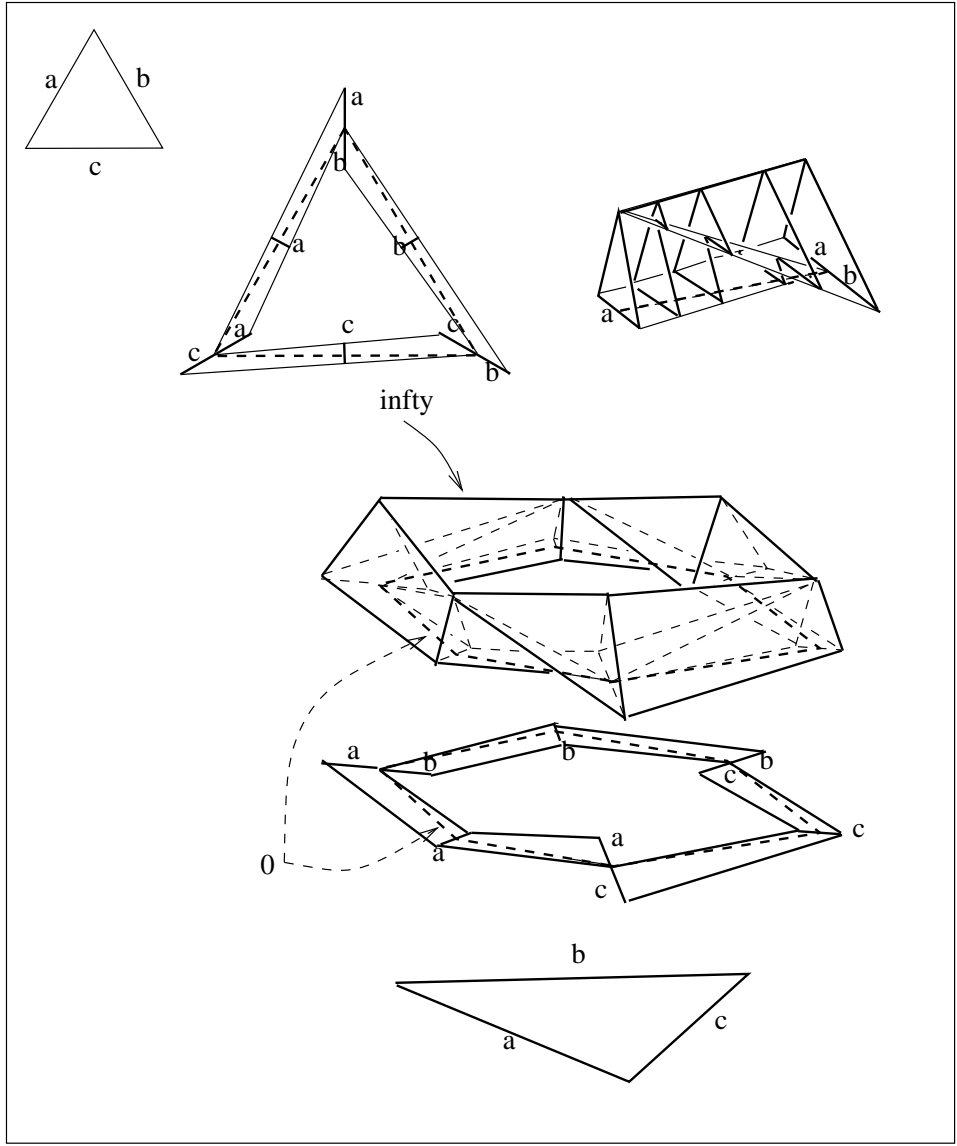}$$
\caption{\label{pic4} The cellular $(S^1,0,\infty)$ model of the
tangent bundle of a triangle }
\end{figure}

To verify 1 and overall naturalness of $\mb G^\infty$ it is useful
to consider the case when $\ms M^n$ is the combinatorial sphere $\ms
S^n$ (see Figure \ref{pic3}).
Let $s$ be a simplex of $\ms S^n$. Then one can imagine $\mb
G^\infty (s)\in \bg R_n$ as follows: glue together all the simplices
of $\ms S^n$ which do not contain $s$. This will be our new ball
of $\mb G^\infty$ marked by $\infty$. The ball is naturally attached
by $\Lk s$ to $\St s$ and altogether they form the sphere $\mb
G^\infty(s)$ with a marked $\infty$-ball. We should mention that
 while $S^n$
is a combinatorial manifold, the ``tangent sphere'' $\mb
G^\infty(s)$ is an abstract manifold, since the $\infty$-ball is
usually non-simplicial. Let $s_0 \subset s_1$ be a pair of simplices
of $\ms S^n$. Then $\St s_0 \supset \St s_1$ When we pass from $\mb
G^\infty (s_0)$ to $\mb G^\infty (s_1)$ the simplices from $\St s_0
\setminus \St s_1$ are dissolve in $\infty$-ball of $\mb
G^\infty (s_1)$. This operation is exactly the morphism $\mb G^\infty_{s_0
\subset s_1}$ and as we see this is exactly an aggregation morphism
from $\bg R_n$. In the  Figure \ref{pic4} we show
the cellular $(S^1,0,\infty)$ model of the tangent bundle of a
triangle obtained by our recipe.
\subsection{$B\bg R_n\approx B\PL_n$}
Now we will discuss the problem of proving the statement 4 of the
Theorem \ref{th1}. Let $X$ be a compact $\PL$ manifold. Consider the
category $\bg R(X)$ of all abstract manifolds $\ms M$ such that the
polyhedron $B\ms M$ is $\PL$ homeomorphic  to $X$. The morphisms of
$\bg R(X)$ are the aggregation morphisms. The Theorem A in
\cite{Mnev:2006} states that
\begin{equation}
\label{e1} B\bg R(X)\approx B\PL(X)
\end{equation}
where $\PL(X)$ is a simplicial group of $\PL$ homeomorphisms of $X$.
The proof of the statement $B\bg R_n\approx B\PL_n$ is a cosmetic
variation of the general scheme developed for (\ref{e1}).

Let $L$ be a $\PL$ polyhedron. We call an \emph{$\bg R(X)$-coloring} of
$L$ the following object: a linear triangulation $\ms K$ of $L$,
$|\ms K|=L$ and an assignment to any vertex of $\ms K$ an abstract
manifold from $\bg R(X)$ and to any 1-simplex of $\ms K$ an
aggregation morphism, in such
a way that all 2-simplices of $\ms K$ become
commutative triangles in $\bg R(X)$. So, the $\bg R(X)$-coloring of
$L$ is just a commutative diagram in $\bg R(X)$ drawn on 2-skeleton
of some triangulation of $ L$. The concordance of two $\bg
R(X)$-colorings $\xi_0,\xi_1$ of $L$ is a coloring of the polyhedron
$L \times [0,1]$, which induces the coloring $\xi_i$ on the $i$-th
side. By  abstract nonsense to prove (\ref{e1}) is the same
as to establish functorial one-to-one correspondence between
isomorphism classes of $\PL$ fiber-bundles on $L$ with  fiber $X$
and concordance classes of $\bg R(X)$ colorings of $L$.

We will menton how we would like but cannot  establish such a
correspondence. This speculation is borrowed  from
\cite{steinberger:1986}. To any $\bg R(X)$-coloring of $L$ we can
apply the construction $\mr{Hocolim}$, and its geometric realization will
produce a triangulated fiber bundle with fiber $X$. On the other
side one can triangulate any fiber bundle with the base $L$. To any
triangulated fiber bundle $ J$ with  base $L$ and fiber $X$ one can
canonically associate (\cite{hatcher:1975},\cite{steinberger:1986})
some $\bg R(X)$-coloring of the first barycentric subdivision of the
base of $ J$. We call this construction by $\mr{Hocolim}^{-1}$. The
composition $\mr {Hocolim} \circ \mr{Hocolim}^{-1}$ applied to a
bundle produces an isomorphic bundle. We would prove (\ref{e1}) in a nice and
short way if we could establish some canonical concordance between
any $\bg R(X)$-coloring $\xi$ of $L$ and the coloring
$\mr{Hocolim}^{-1}\mr{Hocolm}\,\xi$. This would magically eliminate
geometry. Unfortunately there is no way to see such a canonical
concordance. This is the cause of some published and many
unpublished mistakes. From the theory developed in \cite{Mnev:2006}
it follows that $\xi$ and $\mr{Hocolim}^{-1}\mr{Hocolm}\,\xi$ are
concordant, but the concordance is transcendental.

\begin{figure}\caption{\label{polycone00}}
$$\includegraphics[scale=0.7]{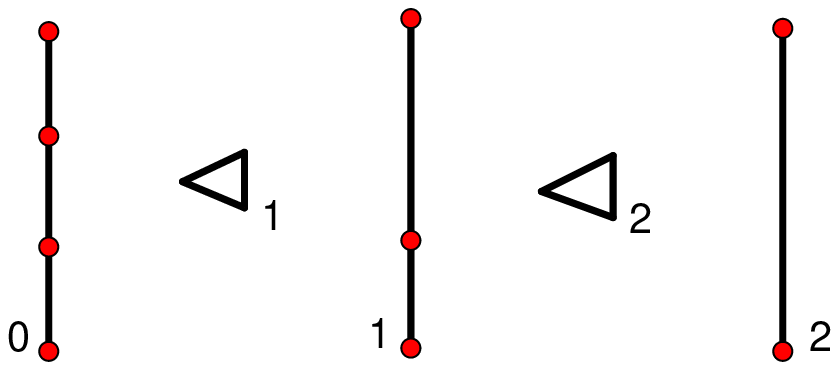}$$
\end{figure}

\begin{figure}\caption{\label{polycone11}}
$$\includegraphics[scale=0.7]{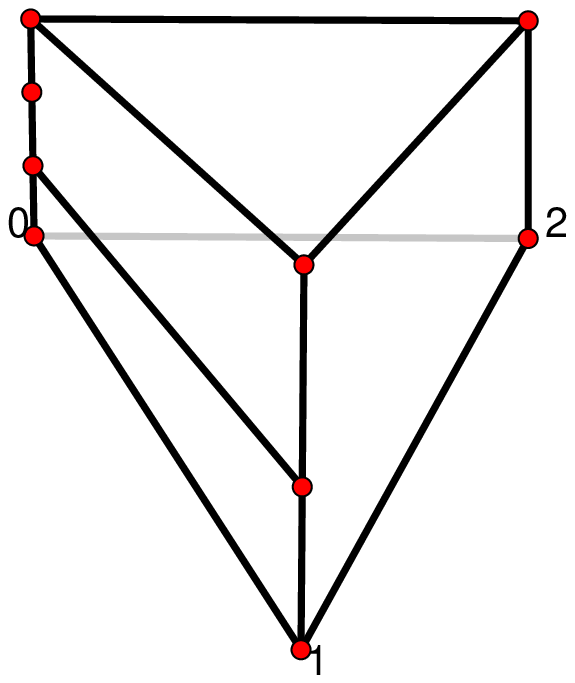}$$
\end{figure}

Instead of using $\mr{Hocolim}$-construction we are constructing a bundle on
$L$ from $\bg R(X)$-coloring of $L$ using traditional construction of
trivializations and  structure homeomorphisms. Let $K$ be a $\bg R
(X)$-colored simplicial complex, $ |K|=L$. The coloring induces
coloring of a $k$-simplex of $K$ by the chain
$$\ms Q_0\rightsquigarrow \ms Q_1\rightsquigarrow...\rightsquigarrow \ms Q_k$$
of abstract aggregations. According to speculations in \S\ref{aggr}
on page \pageref{aggr} one can realize this chain  by the chain
\[ Q=( Q_0 \trianglelefteq  Q_1 \trianglelefteq ... \trianglelefteq Q_k)\]
of geometric aggregations of geometric $\PL$-ball complexes. With
the  chain $Q$ one can associate a ball decomposition
 of the trivial bundle
$X\times \Delta^k\xar{\pi}\Delta^k$ into the  horizontal
\emph{``prisms'' }which are the trivial subbundles with a ball as a
fiber. The drawings \ref{polycone00} and \ref{polycone11}
illustrate the construction of prismatic decomposition on $\pi$ by
the chain of geometric aggregations.
The combinatorics of the coloring associates to any pair of
simplices $s_0 \subset s_1$ in $K$ a combinatorial isomorphism of
two prismatic structures on trivial bundle over $s$. By  Alexander's
trick on can represent all these combinatorial isomorphisms  by
fiberwise structure $\PL$-homeomorphisms of  fiber bundle with the
base $L$ and the fiber $X$. All these structure homeomorphisms map
prisms to prisms. As a result we obtained from $\bg R(X)$-colorings the
class of fiber bundles with unusual structure homeomorphisms -- the
``prismatc'' ones.

In this setup the inverse problem is to learn how to deform the
structure homeomorphisms of arbitrary $\PL$ fiber bundle into the
``prismatic'' form and present a consistent coloring in a
controllable way.

At this point it is useful to recall the proof of the Lemma on
fragmentation of isotopy. This Lemma was proved by Hudson  \cite{Hu}
in the $\PL$ case. It states that for any covering  $U =\{ U_i \}_i$
of a manifold $X$ by open balls and for any $\PL$-homeomorphism $X
\xar{f}X$ which is isotopic to identity there exist a finite
decomposition $f=f_1\circ....\circ f_m$ such that $\forall i \exists
j: \supp f_i \subset U_j$. The proof of the fragmentation lemma
contains more information than the statement. In the proof we pick
arbitrary $\PL$-isotopy $F$ connecting $f$ and identity. Then we
deform $F$ in the class of isotopies with fixed ends to the isotopy
$F'$ of a special form. The isotopy $F'$ corresponds to the chain of
isotopies which are fixed on complements of open balls from $U$.

\begin{figure}\caption{\label{pic5}}
\[\begin{picture}(0,0)%
\includegraphics{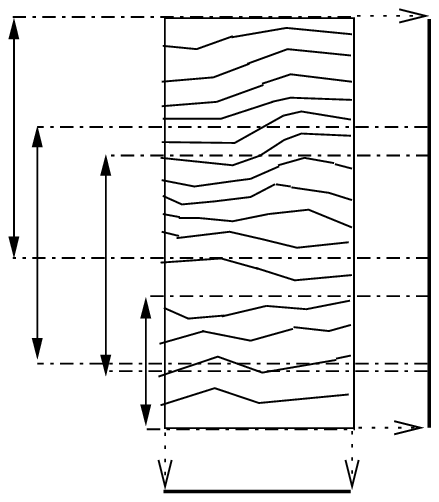}%
\end{picture}%
\setlength{\unitlength}{4144sp}%
\begingroup\makeatletter\ifx\SetFigFont\undefined%
\gdef\SetFigFont#1#2#3#4#5{%
  \reset@font\fontsize{#1}{#2pt}%
  \fontfamily{#3}\fontseries{#4}\fontshape{#5}%
  \selectfont}%
\fi\endgroup%
\begin{picture}(2212,2387)(2094,-15793)
\put(2474,-14339){\makebox(0,0)[lb]{\smash{\SetFigFont{6}{8.4}{\rmdefault}{\mddefault}{\updefault}{\color[rgb]{0,0,0}$U_2$}%
}}}
\put(2729,-15244){\makebox(0,0)[lb]{\smash{\SetFigFont{6}{8.4}{\rmdefault}{\mddefault}{\updefault}{\color[rgb]{0,0,0}$U_4$}%
}}}
\put(2570,-14864){\makebox(0,0)[lb]{\smash{\SetFigFont{6}{8.4}{\rmdefault}{\mddefault}{\updefault}{\color[rgb]{0,0,0}$U_3$}%
}}}
\put(2094,-13978){\makebox(0,0)[lb]{\smash{\SetFigFont{6}{8.4}{\rmdefault}{\mddefault}{\updefault}{\color[rgb]{0,0,0}$U_1$}%
}}}
\put(2895,-15683){\makebox(0,0)[lb]{\smash{\SetFigFont{10}{14.4}{\familydefault}{\mddefault}{\updefault}{\color[rgb]{0,0,0}0}%
}}}
\put(3233,-15793){\makebox(0,0)[lb]{\smash{\SetFigFont{6}{8.4}{\rmdefault}{\mddefault}{\updefault}{\color[rgb]{0,0,0}$I=[0,1]$}%
}}}
\put(4306,-14513){\makebox(0,0)[lb]{\smash{\SetFigFont{6}{8.4}{\rmdefault}{\mddefault}{\updefault}{\color[rgb]{0,0,0}$X$}%
}}}
\put(3923,-15683){\makebox(0,0)[lb]{\smash{\SetFigFont{10}{14.4}{\familydefault}{\mddefault}{\updefault}{\color[rgb]{0,0,0}1}%
}}}
\end{picture}\]
\end{figure}

\begin{figure}\caption{\label{pic6}}
$$\begin{picture}(0,0)%
\includegraphics{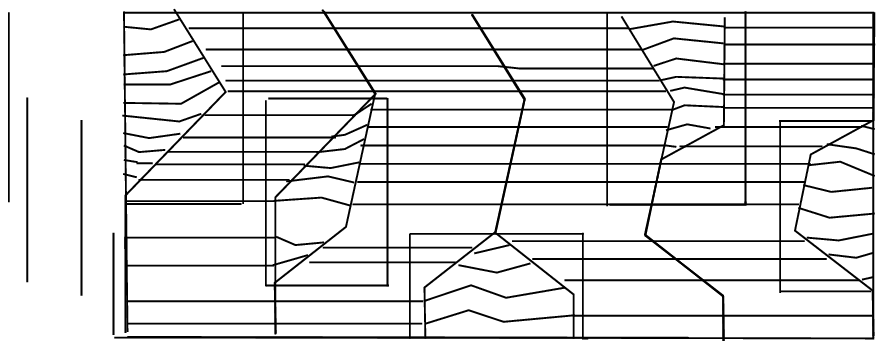}%
\end{picture}%
\setlength{\unitlength}{4144sp}%
\begingroup\makeatletter\ifx\SetFigFont\undefined%
\gdef\SetFigFont#1#2#3#4#5{%
  \reset@font\fontsize{#1}{#2pt}%
  \fontfamily{#3}\fontseries{#4}\fontshape{#5}%
  \selectfont}%
\fi\endgroup%
\begin{picture}(4152,1540)(4231,-16379)
\put(4485,-15598){\makebox(0,0)[lb]{\smash{\SetFigFont{6}{7.2}{\familydefault}{\mddefault}{\updefault}{\color[rgb]{0,0,0}$U_2$}%
}}}
\put(4567,-15990){\makebox(0,0)[lb]{\smash{\SetFigFont{6}{7.2}{\familydefault}{\mddefault}{\updefault}{\color[rgb]{0,0,0}$U_3$}%
}}}
\put(4710,-16290){\makebox(0,0)[lb]{\smash{\SetFigFont{6}{7.2}{\familydefault}{\mddefault}{\updefault}{\color[rgb]{0,0,0}$U_4$}%
}}}
\put(4231,-15376){\makebox(0,0)[lb]{\smash{\SetFigFont{6}{7.2}{\familydefault}{\mddefault}{\updefault}{\color[rgb]{0,0,0}$U_1$}%
}}}
\end{picture}$$
\end{figure}

The isotopy  $F$ is the same thing as a fiberwise homeomorphism
$$X\times {[0,1]} \xar{F} X\times {[0,1]}$$
commuting with the projection on $[0,1]$ and
such that $F_0 = \id$ and $F_1=f$. The homeomorphism $F$ is the same
thing as a one-dimensional foliation $\mc F$ on $X\times [0,1]$
transversal to the fibers of the projection.
(Figure \ref{pic5}).
The homeomorphism $F'$ corresponds to foliation $\mc F'$ with
following property: for any point $b\in [0,1]$ all the points  $x
\in X$ such that the leaf of $\mc F$ ``is not horizontal'' at
$(x,b)$ are contained in an element of $ U$ (Figure \ref{pic6}).
Inspecting  the drawing of  $\mc F'$ one can see that it is possible
to subdivide the base $[0,1]$ into the intervals $u_1,...,u_m$ and
introduce prismatic structure on all subbundles $X\times u_i
\xar{\pi_2} u_i$ such that induced homeomorphisms $F'\lfloor_{u_i}$
are prismatic. So, the construction of the  fragmentation lemma
allows us to deform a fiberwise homeomorphism of the trivial bundle
over interval into the system of prismatic homeomorphisms over the
subdivision of interval. The deformation  $F\rightsquigarrow F'$ has
a canonical form and possesses a coordinate generalization to the
homeomorphisms of the trivial bundle over cube.
So, realizing the program of such a generalization together with the
development of appropriate surgery for fiberwise homeomorphisms
consumes about 100 pages in \cite{Mnev:2006}.

\def\cprime{$'$} \def\cprime{$'$}

\end{document}